\documentclass{amsart}
\usepackage{amsfonts,amssymb,amscd,amsmath,enumerate,verbatim,calc}
\usepackage[all]{xy}

\newcommand{\CM}{Cohen-Macaulay}

\newcommand{\F}{\mathbb{F} }

\newcommand{\yb}{\mathbf{y} }

\newcommand{\q}{\mathfrak{q} }

\newcommand{\C}{\mathcal{C} }
\newcommand{\Z}{\mathbb{Z} }
\newcommand{\Pc}{\mathcal{P} }
\newcommand{\Qc}{\mathcal{Q} }

\newcommand{\db}{\mathbf{d} }
\newcommand{\xb}{\mathbf{x} }
\newcommand{\rt}{\rightarrow}

\newcommand{\G}{\mathbb{G}^\bullet }
\newcommand{\Sq}{\mathcal{S}{q}}
\newcommand{\fP}{\mathfrak{P}}
\newcommand{\fp}{\mathfrak{p}}

\newcommand{\wh}{\widehat }

\newcommand{\image}{\operatorname{image}}
\newcommand{\coker}{\operatorname{coker}}

\newcommand{\Supp}{\operatorname{Supp}}
\newcommand{\grade}{\operatorname{grade}}
\newcommand{\depth}{\operatorname{depth}}
\newcommand{\height}{\operatorname{height}}
\newcommand{\ann}{\operatorname{ann}}
\newcommand{\Ass}{\operatorname{Ass}}

\newcommand{\Spec}{\operatorname{Spec}}

\newcommand{\Hom}{\operatorname{Hom}}
\newcommand{\cHom}{\operatorname{Hom}}

\newcommand{\Ext}{\operatorname{Ext}}

\theoremstyle{plain}

\newtheorem{theorem}{Theorem}[section]
\newtheorem{corollary}[theorem]{Corollary}
\newtheorem{lemma}[theorem]{Lemma}
\newtheorem{proposition}[theorem]{Proposition}

\theoremstyle{definition}
\newtheorem{definition}[theorem]{Definition}
\newtheorem{s}[theorem]{}
\newtheorem{remark}[theorem]{Remark}
\newtheorem{example}[theorem]{Example}

\theoremstyle{remark}

\begin{document}

\title[Steenrod algebras]{Local cohomology over Steenrod algebras}
\author{Tony~J.~Puthenpurakal}
\date{\today}
\address{Department of Mathematics, Indian Institute of Technology Bombay, Powai, Mumbai 400 076, India}

\email{tputhen@gmail.com}

 \begin{abstract}
Let $\F_q$ be a finite field with $q = p^s$ elements. Let $S = \bigoplus_{n \geq 0} S_n$ be a non-negatively  graded commutative,  Noetherian $\F_q$-algebra. Assume $S$ has a structure of an unstable Steenrod algebra. Let $M$ be a finitely generated graded unstable $S$-module. Let $I$ be a homogeneous ideal in $S$. We prove that for $i \geq 0$ the local cohomology modules $H^i_I(M)$ has a structure of a Cartan $S$-module. As application we give a direct and  considerably simpler proof of D. Bourguiba and S. Zarati's generalization of Landweber-Stong conjecture. We also give a necessary and sufficient condition for $S$ to have $S_r$-property of Serre in terms of Dickson  polynomials.
\end{abstract}
 \maketitle
\section{introduction}
Let $\F_q$ be a finite field with $q = p^s$ elements. Let $S = \bigoplus_{n \geq 0} S_n$ be a non-negatively  graded  Noetherian $\F_q$-algebra with $S_0 = \F_q$. Assume $S$ has a structure of an unstable Steenrod algebra. Recall that this means there exists $\F_q$ linear maps $\Pc^i \colon S \rt S$ with $\Pc^0$ to be the identity and satisfy the Cartan identity:
$$\Pc^k(uv) = \sum_{i+j = k}\Pc^i(u)\Pc^j(v). $$
We  have for homogeneous $u$, $\Pc^{\deg u}(u) = u^q$ and $\Pc^{i}(u) =  0$ for $i > \deg u$.
We also have that $\Pc^r(S_n) \subseteq S_{n + r(q-1)}$ for all $n \in \Z$.

A large class of examples of unstable Steenrod algebras are rings of invariants $S$ of a finite group $G \subseteq GL_d(\F_q)$.  We note that $GL_d(\F_q)$ is a finite group  and so $S$ contains $D^*(d)$ the ring of invariants of $GL_d(\F_q)$. Dickson proved that $D^*(d)$ is a polynomial ring with
$D^*(d) = \F_q[ \db_{d,0}, \db_{d, 1}, \ldots, \db_{d, d-1}]$ with $\deg \db_{d,i} = q^d - q^i$. We follow the exposition given in \cite{SL}.

 One of the landmark results in modular invariant theory is the proof  by D. Bourguiba and S. Zarati, see \cite{bz}, of
the Landweber-Stong conjecture \cite{LS}. This result states that $\depth S \geq r$ if and only if $\db_{d,d-1}, \cdots, \db_{d, d-r}$ is an $S$-regular sequence. In our paper \cite{P} we gave an  alternate and considerably simpler  proof of this conjecture.

D. Bourguiba and S. Zarati prove a considerably more general result than
the Landweber-Stong conjecture, see \cite[1.3]{bz} A natural question was whether our techniques generalize. In this paper we show that our techniques from \cite{P} generalize and give an alternate and considerably simpler proof of D. Bourguiba and S. Zarati's result.

The crucial point is the notion of a Cartan $S$-module which was introduced in \cite{P}.
 \begin{definition}
By a Cartan $S$-module we mean a graded  $S$-module $M$ equipped with $\F_q$-linear maps $\Qc^i_M \colon M \rt M$ for $i \geq 0$ with $\Qc^0_M $ being the identity such that for any $s \in S$ and $m \in M$ we have for all $r \geq 1$
\begin{enumerate}
  \item $$ \Qc^r_M(sm) = \sum_{i+j = r}\Pc^i(s)\Qc^j_M(m).$$
  \item  $\Qc^r_M(M_n) \subseteq M_{n + r(q-1)} $ for all $n \in \Z$.
\end{enumerate}
\end{definition}

 It is easily seen that if $M$ is a Cartan $S$-module then $\ann_S M$ is a $\Pc^*$-invariant ideal of $S$, see \ref{ann-p*}.
 The main result which enables all our applications is the following:
 \begin{theorem}
 \label{jim} Let $S = \bigoplus_{n \geq 0}S_n$ with $S_0 = \F_q$ be a Noetherian $\F_q$ graded algebra. Assume $S$ has a structure of an unstable Steenrod algebra. Let $M$ be a Cartan $S$-module (not necessarily finitely generated). Let $I$ be a homogeneous ideal in $S$ and let $\xb = x_1, \ldots, x_m$ be a set of homogeneous generators of $I$. Then for $i \geq 0$ the local cohomology modules $H^i_I(M)$ have the structure of a Cartan $S$-module (depending on $\xb$).
 \end{theorem}

\begin{remark}\label{lannes}
Let $S$ be an unstable algebra over the Steenrod algebra. Let $\dim S = d$. Then there exists an embedding  of unstable Steenrod algebras $$\psi \colon \F_q[ \db_{d,0}^{q^t}, \db_{d, 1}^{q^t}, \ldots, \db_{d, d-1}^{q^t}]\rt S$$ such that $S$ is a finite $\F_q[ \db_{d,0}^{q^t}, \db_{d, 1}^{q^t}, \ldots, \db_{d, d-1}^{q^t}]$-module via $\psi$. For a proof of this result see \cite[1.2]{bz}.
We will set  $\F_q[ \db_{d,0}^{q^t}, \db_{d, 1}^{q^t}, \ldots, \db_{d, d-1}^{q^t}]$ as $D(d)^{q^t}$. We will fix such an embedding for $S$. By abuse of notation we consider the embedding  as an \emph{inclusion}.
\end{remark}

\emph{Applications:}

I. We give a considerably simpler and a direct proof of  Bourguiba and  Zarati's result.

\begin{theorem}
\label{bz-gen}(with setup as in \ref{lannes}) Let $M$ be a finitely generated Cartan $S$-module. The following assertions are equivalent
\begin{enumerate}[\rm (i)]
  \item $\depth M \geq r$.
  \item $\db_{d,d-1}^{q^t}, \cdots, \db_{d, d-r}^{q^t}$ is an $M$-regular sequence.
\end{enumerate}
\end{theorem}

 II.  \emph{Associate primes of Cartan modules}\\
An ideal $I$ in $S$ is said to   be $\Pc^*$-invariant if $\Pc^*(I) \subseteq I$ for all $i \geq 0$. We prove
\begin{theorem}
\label{ass}(with setup as in \ref{lannes}) Let $M$ be a finitely generated Cartan $S$-module. If $P \in \Ass M$ then $P$ is a $\Pc^*$ invariant.
\end{theorem}
Proof of Theorem \ref{ass} requires that if $M$ is a Cartan $S$-module then for any ideal $J$ of $S$ we have $H^0_J(M)$ is a Cartan submodule of $M$. This follows from Theorem \ref{jim}. However we also give a direct proof. Theorem \ref{ass} is \emph{not} true if $M$ is not finitely generated. We give a counter-example, see \ref{non-ass}.

III. \emph{$S_r$ property}: \\
Let $T$ be a Noetherian ring $T$ and let $M$ be a finitely generated $T$-module. Set $\depth 0 = + \infty$. Recall $M$ satisfies $S_r$ property of Serre if
$$ \depth M_P \geq \min\{\dim M_P, r \} \quad \text{for all primes $P$ of $T$}. $$

A natural question is what happens if $\db_{d,0}^{q^t}, \db_{d, 1}^{q^t},\cdots \db_{d,r-1}^{q^t}$ is an $M$-regular sequence. We extend our result in \cite[1.6]{P} to prove
\begin{theorem}
\label{s_r} (with setup as in \ref{lannes}) Let $M$ be a finitely generated Cartan $S$-module. Assume $\Supp(M) = \Spec(S)$  and that $S$ is equi-dimensional. The following assertions are equivalent:
\begin{enumerate}[\rm (i)]
  \item $M$ satisfies $S_r$-property of Serre.
  \item  $\db_{d,0}^{q^t}, \db_{d, 1}^{q^t},\cdots \db_{d,r-1}^{q^t}$ is an $M$-regular sequence.
\end{enumerate}
\end{theorem}

IV. \emph{Cohomological annihilators:}
\s \label{roti} If $M$ is a finitely generated Cartan $S$-module then  $H^i_{S_+}(M)$ is a Cartan $S$-module. It is not difficult to prove that $\height \ann H^i_{S_+}(M) \geq d -i$, see \ref{prev}. As $\ann H^i_{S_+}(M)$ is $\Pc^*$-invariant it follows that $\db_{d,0}^{q^t} \in \sqrt{\ann H^i_{S_+}(M)} $ for $i < d = \dim S$. Thus some power $\db_{d,0}^{q^t}$ works as a cohomological annihilator. We show
\begin{corollary}
\label{rob} (with hypotheses as in \ref{lannes}). Let $\G = 0 \rt G^0 \rt G^1 \rt \ldots \rt G^r \rt 0$ be a complex of free graded $S$-modules and homogeneous maps such that $H^i(\G)$ has finite length for all $i$. Let $M$ be a finitely generated Cartan $S$-module. Then $\db_{d,0}^{q^t}\in \sqrt{\ann_S H^i(\G\otimes M)}$ for $0 \leq i \leq d-1$.
\end{corollary}
\begin{proof}
  This follows from an exercise problem \cite[8.1.6]{BH} and the discussion in \ref{roti}.
\end{proof}
In \cite[11.5.5]{SL} it is proved that if $S$ is a ring of invariants then $S_{\db_{d_0}}$ is \CM. We extend this result.
\begin{corollary}
  \label{cm} (with hypotheses as in \ref{lannes}).     Let $M$ be a finitely generated Cartan $S$-module. Then $M_{\db_{d,0}^{q^t}}$ is \CM.
\end{corollary}
\begin{remark}\label{sbi}
  If $\dim M < d$ and $M$ is a Cartan module then by \ref{bella}, $\db_{d_0}^{q^t} \in \sqrt{\ann M}$. So $M_{\db_{d,0}^{q^t}} = 0$. Thus Corollary \ref{cm} has content only if
   $\dim M =  \dim S = d$.
\end{remark}

V. \emph{Monomial ideals:}\\
Let $R = \F_q[X_1,\ldots, X_d]$ with the natural Steenrod algebra structure. Let $I$ be a monomial ideal in $R$. It is not difficult to show that $I$ is a $\Pc^*$-invariant ideal of $R$, see \ref{s-mono}. As an application of this observation we prove:
\begin{theorem}
\label{mono} Let $k$ be a field of characteristic $p > 0$ containing $\F_q$. Set $R = \F_q[X_1,\ldots, X_d]$ and $T = k[X_1,\ldots, X_d]$. Let $I$ be a monomial ideal in $T$ and let $J$ be a homogeneous ideal in $R$. Assume $M = H^i_{JT}(T/I)$ is  finitely generated as a $T$-module for some $i \geq 0$. If $P \in \Ass_T M$ then $P$ is generated by linear polynomials in $T$.
\end{theorem}

VI. \emph{Finite generation of cohomology:}\\
Let $S$ be an unstable Steenrod algebra. If $\dim S \geq 2$ then there are  infinitely many homogeneous prime ideals. Note there exists only finitely many  $\Pc^*$-invariant prime ideals.
For a prime ideal $P$ of $S$ define
$$f_P = \min \{ i \mid H^i_P(S) \ \text{is not finitely generated $S$-module} \}. $$
We prove
\begin{theorem}
\label{bound}(with setup as above). Assume $P$ is not a $\Pc^*$-invariant prime ideal. Then $f_P \leq \depth S_P$.
\end{theorem}
\begin{remark}
For prime ideals which are  $\Pc^*$-invariant, the bound in Theorem \ref{bound} need not hold. See \ref{counter-bound} for an example.
\end{remark}

We now describe in brief the contents of this paper.  In section two we discuss a few preliminary results that we need. In Section three we discuss our notion of Cartan modules and discuss a few preliminary properties. In section four we give an alternate description of Cartan $S$-modules which is convenient for us.
In the next section
 discuss a Cartan $S$-module structure on $M_x$ where $x \in S$ is homogeneous and $M$ is a Cartan $S$-module.  In section six we 
    prove Theorem \ref{jim}. In the next section we give a proof of Theorem \ref{bz-gen}. In section eight we give proof of Theorem \ref{ass}. In the next section we give proof of our applications III. Finally in section ten we give a proof of Corollary \ref{cm}. In the next section we give a proof of Theorem \ref{mono}. Finally in section twelve we give a proof of Theorem \ref{bound}.

\section{preliminaries}
\textbf{I:} \emph{Steenrod operators:}\\
We give the following well-known construction of Steenrod operators on polynomial rings and ring of invariants.
We follow the exposition given in \cite{SL}.
\s Let $p$ be a prime, $q = p^s$ and let $V$ be a $\F_q$-vector space of dimension $d$.
Define $P(\xi)\colon \F_q[V] \rt \F_q[V][\xi]$ (with $\deg \xi = 1-q$) by the rules:
\begin{enumerate}
  \item $P(\xi)(v)  = v + v^q\xi$  for all $v \in  V^*$.
  \item $P(\xi)(u + w) = P(\xi)(u) + P(\xi)(w)$ for all $u,w \in F_q[V]$.
  \item $P(\xi)(uw) = P(\xi)(u)P(\xi)(w)$ for all $u,w \in F_q[V]$.
  \item $P(\xi)(1) = 1$.
   \end{enumerate}
  We note that $P(\xi)$ is a ring homomorphism of degree zero.

By separating out homogeneous components we obtain $\F_q$ linear maps \\ $\Pc^i \colon \F_q[V] \rt \F_q[V]$ by the requirement
$$P(\xi)(f) = \sum_{i \geq 0}\Pc^i(f)\xi^i.$$
Note $\Pc^0$ is the identity.
The operations $\Pc^i$ are called the \emph{Steenrod reduced power operations} over $\F_q$ and when $p = 2$, also denoted by $\Sq^i$ and referred to as \emph{Steenrod squaring operations}.
We will use the following Cartan formulae:
$$\Pc^k(uv) = \sum_{i+j = k}\Pc^i(u)\Pc^j(v). $$

We  have for homogeneous $u$, $\Pc^{\deg u}(u) = u^q$ and $\Pc^{i}(u) =  0$ for $i > \deg u$.
We also have that $\Pc^r(\F_q[V]_n) \subseteq \F_q[V]_{n + r(q-1)}$ for all $n \in \Z$.

\s Let $G \subseteq GL(V)$ be a finite group. Let $S = \F_q[V]^G$. We define an action of $G$ on $\F_q[V][\xi]$ by giving the trivial action on $\xi$. We note that $\sigma(P(\xi)(r)) = P(\xi)(\sigma(r))$ for all $\sigma \in G$. Thus the maps $\Pc^i$ restrict to maps from $S$ to $S$. So $S$ is an unstable Steenrod algebra.

\begin{proposition}\label{s-mono}
Let $R = \F_q[X_1, \ldots, X_d]$ be the polynomial ring in $d$-variables over $\F_q$ and give it the unstable Steenrod algebra structure as defined above. Let $I$ be a monomial ideal in $R$. Then
\begin{enumerate}[\rm (1)]
  \item $I$ is a $\Pc^*$-invariant ideal.
  \item $R/I$ is an unstable Steenrod algebra.
\end{enumerate}
\end{proposition}
We will need the following easily proved result.
\begin{lemma}\label{sum-prd}
Let $S$ be an unstable Steenrod algebra and let $I, J$ be $\Pc^*$-invariant ideal. Then $I + J$ and $IJ$ are $\Pc^*$-invariant ideals. \qed
\end{lemma}
We now give
\begin{proof}[Proof of Proposition \ref{s-mono}]
(1) By using \ref{sum-prd} repeatedly it suffices to prove that $(X_i)$ is  $\Pc^*$-invariant. But this is obvious.

(2) This follows from (1).
\end{proof}
\section{Cartan Modules}
In this section we describe our notion of Cartan modules. We work with the setup as in \ref{lannes}. Let $\Pc^i \colon S \rt S$ be the $i^{th}$-Steenrod operator on $S$.

\begin{definition}
By a Cartan $S$-module we mean a graded $S$-module $M$ equipped with $\F_q$-linear maps $\Qc^i_M \colon M \rt M$ for $i \geq 0$ with $\Qc^0_M $ being the identity such that for any $s \in S$ and $m \in M$ we have for all $r \geq 1$
\begin{enumerate}
  \item $$ \Qc^r_M(sm) = \sum_{i+j = r}\Pc^i(s)\Qc^j_M(m).$$
  \item $\Qc^r_M(M_n) \subseteq M_{n + r(q-1)}$ for all $n \in \Z$.
\end{enumerate}
\end{definition}
Next we define morphisms of Cartan $S$-modules.
\begin{definition}\label{def-cart-morph}
  Let $M, N$ be Cartan $S$-modules. By a morphism $f \colon M \rt N$ of Cartan $S$-modules we mean a homogeneous $S$-linear map $f \colon M \rt N$ such that the following diagram commutes for all $i \geq 0$.

\[
  \xymatrix
{
M
\ar@{->}[r]^{f}
\ar@{->}[d]^{\Qc_M^i}
 & N
\ar@{->}[d]^{\Qc_N^i}
\\
M
\ar@{->}[r]^{f}
 & N
\
 }
\]
\end{definition}
\s Let $M, N$ be Cartan $S$-modules. We denote
\[
\cHom_{\C(S)}(M, N) = \{ f \colon M \rt N \mid f \ \text{is a morphism of Cartan $S$-modules}\}.
\]
It is evident that $\cHom_{\C(S)}(M, N)$ is an abelian group.
The following result follows easily from the definition of morphism of Cartan $S$-modules.
\begin{proposition}
\label{kic} Let $M, N$ be Cartan $S$-modules and let $f \colon M \rt N$ be a morphism of Cartan $S$-modules. Then $\ker f, \image f, \coker f$ are Cartan $S$-modules. Furthermore the natural inclusions $\ker f \rt M$, $\image f \rt N$ and the surjection $N \rt \coker f$ are morphisms of Cartan $S$-modules. \qed
\end{proposition}

We will need a few preliminaries.
\begin{proposition}\label{add}
Let $M, N$ be Cartan $S$-modules. Then $M \oplus N$  also has a natural Cartan $S$-module structure such that
\begin{enumerate}[\rm (1)]
  \item the natural maps $M \rt M\oplus N$, $N \rt M\oplus N$, $M\oplus N \rt M$ and $M\oplus N \rt N$ are morphisms of Cartan $S$-modules.
  \item if $U$ is a Cartan $S$-module and $f\colon M \rt U$ and $g \colon N\rt U$ are morphism of Cartan $S$-modules then $f\oplus g \colon M\oplus N \rt U$ is a morphism of Cartan $S$-modules.
  \item if $U$ is a Cartan $S$-module  and $\alpha \colon U \rt M$ and $\beta \colon U \rt N$ are morphism of Cartan $S$-modules then $(\alpha, \beta) \colon  U \rt M\oplus N$ is a morphism Cartan $S$-modules.
\end{enumerate}
\end{proposition}
\begin{proof}
Define $\Qc^r_{M\oplus N} \colon M\oplus N \rt M\oplus N$ to be
$$\Qc^r_{M\oplus N}(m, n) = (\Qc^r(m), \Qc^r(n)).$$
It is readily verified that this gives a Cartan $S$-module structure on $M\oplus N$.

The assertions (1), (2) and (3) hold trivially.
\end{proof}
We also need the following easily proved result.
\begin{proposition}\label{map-add}
Let $M_i, N_j$ be Cartan $S$-modules for $1 \leq i \leq m$ and $1\leq j \leq n$. Set $M = \bigoplus_{i = 1}^{m}M_i$ and $N = \bigoplus_{j = 1}^{n}N_j$. Give $M, N$ the induced Cartan $S$-module structure. Let $\psi = (\psi_{i,j})\colon M \rt N$ where $\psi_{i,j} \colon M_i \rt N_j$. If each $\psi_{i,j}$ is morphism of Cartan $S$-modules then $\psi$ is a morphism of Cartan $S$-modules. \qed
\end{proposition}
\s Let $\C(S)$ denote the category of Cartan $S$-modules and Cartan morphisms. Then by \ref{kic}, \ref{add}, \ref{map-add} it follows that $\C(S)$ is an abelian category.
As a consequence we obtain
\begin{corollary}
\label{complex-cart} Let $\G $ be a complex of Cartan $S$-modules and morphisms.  Then $H^i(\G)$ are  Cartan $S$-modules for all $i \in \Z$.
\end{corollary}

\s \label{cleo} An ideal $I$ of $S$ is said to be \emph{$\Pc^*$-invariant} if $\Pc^i(I)\subseteq I$ for all $i \geq 0$. We note that if $I$ is $\Pc^*$-invariant then so is $\sqrt{I}$, see \cite[11.2.5]{SL}.

The following Lemma is critical. We proved  it in \cite[3.5]{P} for rings of invariants. Same proof works in general.
\begin{lemma}\label{ann-p*}
Let $M$ be a Cartan $S$-module. Then $\ann M$ is a $\Pc^*$-invariant ideal of $S$.
\end{lemma}

\s \label{dickson} The non-zero $\Pc^*$-prime ideals of $D(d)^{q^t}$  are $(\db_{d,0}^{q^t}, \cdots, \db_{d,i}^{q^t})$ for $i = 0, \ldots, d -1$; see \cite[9.2.7]{NS}.
We need the following result. We proved it for ring of invariants, see \cite[2.8]{P}. The same proof works in general.
\begin{lemma}\label{bella}
(with setup as in \ref{lannes}). Let $J$ be a non-zero  radical $\Pc^*$-invariant ideal of $S$. If $\dim S/J = d - g$ then $\db_{d,0}^{q^t}, \cdots, \db_{d, g -1}^{q^t} \in J$.
\end{lemma}

\section{An alternate description of Cartan modules}
In this section we give an alternate description of Cartan modules  which will be useful later.
Throughout $S$ is an unstable Steenrod $\F_q$-algebra.

\s \emph{The total Steenrod operation}. This is the map defined by $\Pc(\xi) \colon S \rt S[\xi]$ where $\Pc(\xi)(u) = \sum_{i \geq 0}\Pc^i(u)\xi^i$.
We note that $\Pc(\xi)$ is a ring homomorphism. Furthermore if we set $\deg(\xi) = 1-q$ then $\Pc(\xi)$ is a homogeneous ring homomorphism.
We note that $S[\xi]$ has elements of negative degrees (as $\deg(\xi) = 1-q$)  and so is a $\Z$-graded ring.

\s \emph{Important point:}(with setup as above). As $\Pc(\xi) \colon S \rt S[\xi]$ is a ring homomorphism, it induces a \emph{new} $S$ algebra structure on $S[\xi]$.

\begin{definition}
  Let $M$ be an $S$-module and let $N$ be an $S[\xi]$-module and let $f\colon M \rt N$ be an additive homogeneous graded map. We say $f$ is a Steenrod map if $f$ is $S$-linear with $N$ given an $S$-module structure via $\Pc(\xi)$, i.e., $f(u m) = \Pc(\xi)(u)f(m)$ for all $u\in S$ and  for all $m \in M$.
\end{definition}

\s We need to construct a graded subring $S[[\xi]]^*$ of $S[[\xi]]$.  For a homogeneous element $s \in S$ set $|s| =\deg(s)$.
For $n \in \Z$, set
\[
S[[\xi]]^*_n = \{\alpha =  \sum_{i \geq 0} \alpha_i \xi^i \mid  \alpha_i \text{is homogeneous and} \ |\alpha_i| = n + i(q-1)\}.
\]
It is readily verified that $S[[\xi]]^*_n$ is an additive subgroup of $S[[\xi]]$. Furthermore it is easy to check that
$$ S[[\xi]]^*_n S[[\xi]]^*_m \subseteq S[[\xi]]^*_{m+n}.$$
Set $S[[\xi]]^* = \bigoplus_{n \in \Z}S[[\xi]]^*_n $. It is a graded ring. It also contains $S[\xi]$ as a graded subring.

\s Let $M$ be a graded $S$-module (not necessarily finitely generated). We need to construct a graded $S[[\xi]]^*$ submodule $ M[[\xi]]^*$ of $M[[\xi]]$.  For a homogeneous element $m \in M$ set $|m| =\deg(m)$.
For $n \in \Z$, set
\[
M[[\xi]]^*_n = \{\beta =  \sum_{i \geq 0} \alpha_i \xi^i \mid  \beta_i \text{is homogeneous in $M$ and} \ |\beta_i| = n + i(q-1)\}.
\]
It is readily verified that $M[[\xi]]^*_n$ is an additive subgroup of $M[[\xi]]$. Furthermore it is easy to check that
$$ S[[\xi]]^*_n M[[\xi]]^*_m \subseteq M[[\xi]]^*_{m+n}.$$
Set $M[[\xi]]^* = \bigoplus_{n \in \Z}M[\xi]]^*_n $. It is a graded $S[[\xi]]^* $-module. It also contains $M[\xi]$ as  $S[\xi]$-submodule.

We now give our characterization of Cartan modules.
\begin{theorem}
\label{char-cartan}(with setup as above). Give $S[[\xi]]^*$ modules a $S[\xi]$-module structure via the inclusion $S[\xi] \hookrightarrow S[[\xi]]^*$.
Let $M$ be a graded $S$-module. The following assertions are equivalent:
\begin{enumerate}[\rm (a)]
  \item $M$ is a Cartan $S$-module.
  \item There exists a Steenrod map $\Qc_M(\xi)\colon M \rt M[[\xi]]^*$ with $\Qc_M(\xi)(m)_0 = m$ for all $m \in M$.
\end{enumerate}
\end{theorem}
\begin{proof}
(a) $\implies$ (b).\\
Define $\Qc_M(\xi) \colon M \rt M[[\xi]]^*$ by
$$\Qc_M(\xi)(m) = \sum_{n \geq 0}\Qc_M^n(m)\xi^i,$$
with $\deg \xi = 1-q$.
We note that $\Qc_M$ takes homogeneous elements of degree $r$ to degree $r$ elements in $M[[\xi]]^*$.
Let $u \in S$ and $m \in M$ be homogeneous
We have \begin{align*}
         \Pc(\xi)(u)\Qc_M(m)  &= \left(\sum_{n \geq 0}\Pc^n(u)\xi^n \right) \left(\sum_{n \geq 0}\Qc^n_M(m)\xi^n \right) \\
          &= \sum_{n \geq 0}\left(\sum_{i+j = n}\Pc^i(u)\Qc^j_M(m)\right)\xi^n \\
           &= \Qc_M(um), \quad \text{by Cartan's identity}.
        \end{align*}
 The result follows.

  (b) $\implies$ (a).\\
  Set $$\Qc_M(\xi)(m) = \sum_{r \geq 0}\Qc_M^r(m)\xi^r$$
  Then $\Qc_M^r \colon M \rt M$ are $\F_q$-linear maps. Furthermore as $\Qc_M$ is a homogeneous map and $|\xi| = 1 -q$ it follows that $\Qc^r_M(M_n) \subseteq M_{n + r(q-1)} $ for all $n \in \Z$.
  By assumption $\Qc_M^0(m) = m$ for all $m \in M$.
  As $\Qc_M$ is a Steenrod map we have for
   $u \in S$, $m \in M$  homogeneous
   \begin{align*}
     \Qc_M(um) &= \Pc(\xi)(u)\Qc_M(m) \\
     &=  \left(\sum_{n \geq 0}\Pc^n(u)\xi^n \right) \left(\sum_{n \geq 0}\Qc^n_M(m)\xi^n \right) \\
     &= \sum_{n \geq 0}\left(\sum_{i+j = n}\Pc^i(u)\Qc^j_M(m)\right)\xi^n.
   \end{align*}
   The Cartan identity follows.
\end{proof}

\s\label{loc-ring-m} Let $x$ be a non-nilpotent homogeneous element in $S$. Then we can form the $S[[\xi]]^* $- graded module $S_x[[\xi]]^*$. It is readily verified that $S_x[[\xi]]^*$ is a graded ring. Furthermore  the natural
homogeneous map $S[[\xi]]^* \rt S_x[[\xi]]^* $ is a ring homomorphism.
Similarly we can show that if $M$ is a graded $S$-module then  $M_x[[\xi]]^*$ is a graded $S_x[[\xi]]^* $-module.

We need the following result.
\begin{lemma}
\label{inv}Let $S$ be a Steenrod $\F_q$-algebra. Let $x \in S$ be a homogeneous non-nilpotent element. We have $\Pc(\xi)(x) = x + \sum_{i = 1}^{|x|}\Pc^i(x)\xi^i$ is invertible in $S_x[[\xi]]$.
Then \begin{enumerate}[\rm (1)]
       \item $\Pc(\xi)(x)^{-1} \in S_x[[\xi]]^*$.
       \item $\Pc(\xi)(x)^{-1}_0 = 1/x$.
     \end{enumerate}

\end{lemma}
\begin{proof}
(1) and (2).  Let $|x| = c$. We show that $\Pc(\xi)(x)^{-1} \in S_x[[\xi]]^*_{-c}$.
  Set
  $$ \Pc(\xi)(x)^{-1} = \sum_{n \geq 0}\alpha_n\xi^n. $$
  We prove by induction that  $|\alpha_n| = n(q-1) -c$ for all $n \geq 0$.

  We have $ \Pc(\xi)(x)  \Pc(\xi)(x)^{-1}  = 1$. So we have
 $$ \sum_{n \geq 0}\left(\sum_{i+j = n}\Pc^i(x)\alpha_j\right)\xi^n   = 1.$$
 Comparing coefficients we have $x\alpha_0 = 1$. So $\alpha_0 = 1/x$ has degree $-c$. This also verifies (2),

 Let $n \geq 1$. By induction suppose  $|\alpha_r| = r(q-1) -c$ for $r = 0, 1, \ldots, n-1$.
 We have $\sum_{i+j = n}\Pc^i(x)\alpha_j = 0$.
 So
 \[
 x\alpha_n +  \sum_{i+j = n, j\neq n}\Pc^i(x)\alpha_j = 0.
 \]
 We have $$|\Pc^i(x)\alpha_j|  = |\Pc^i(x)||\alpha_j| = i(q-1) + c + j(q-1)-c = n(q-1). $$
 So $|\alpha_n| = n(q-1) -c$. The result follows by induction.
\end{proof}

\section{The Cartan $S$-module structure on $M_x$}
In this section we prove two critical results.
Throughout $S$ is a Steenrod $\F_q$-algebra, $M$ is a Cartan $S$-module and $x \in S$ is homogeneous.
 We first prove
\begin{theorem}
\label{loc}(with setup as in \ref{lannes}).  Then there exists a  Cartan $S$-module structure on $M_x$.
\end{theorem}
The next result is also critical.
\begin{theorem}
\label{maps}(with setup as in \ref{lannes}). Let $x, y$ be homogeneous. Then the natural map $M_x \rt M_{xy}$ is morphism of Cartan $S$-modules (the Cartan module structures on $M_x, M_{xy}$ are as constructed in Theorem \ref{loc}).
\end{theorem}
We first give
\begin{proof}[Proof of Theorem \ref{loc}]
We may assume that $M_x\neq 0$, otherwise there is nothing to prove. In particular $x$ is a non-nilpotent element of $S$.
We note that $M_x[[\xi]]$ is a graded  $S_x[[\xi]]^*$-module and so an $S[\xi]$-module via the natural map $S[\xi] \rt S_x[[\xi]]$.
Let $\Qc_M (\xi)\colon M \rt M[[\xi]]^*$ be the total Steenrod operation on $M$.
We construct a Steenrod map $\Qc_{M_x}(\xi) \colon M_x \rt M_x[[\xi]]^*$ as
$$\Qc_{M_x}(\xi)(\frac{m}{x^n})   =  \Pc(\xi)(x)^{-n}\Qc_M(\xi)(m). $$
We note that the right hand side always belongs to $M_x[[\xi]]^*$, see \ref{inv} and \ref{loc-ring-m}.
We first show that the formula is well defined. Suppose $m_1/x^n = m_2/x^l$. Then there exists $r \geq 0$ such that $x^{l+r}m_1 = x^{n+r}m_2$. As $M$ is a Cartan $S$-module we get
$\Pc(\xi)(x)^{l+r}\Qc_M(m_1) = \Pc(\xi)(x)^{n+r}\Qc_M(m_2)$. So we obtain
$$\Pc(\xi)(x)^{-n}\Qc_M(\xi)(m_1) = \Pc(\xi)(x)^{-l}\Qc_M(\xi)(m_2).$$
Thus we have a well-defined map  $\Qc_{M_x} \colon M_x \rt M_x[[\xi]]^*$.
We note that
$$\Qc_{M_x}(\xi)(\frac{m}{x^n})_0 = \Pc(\xi)(x)^{-n}_0\Qc_M(\xi)(m)_0 = \frac{m}{x^n}$$
using the definition of Cartan modules and \ref{inv}. Also  note that $\Qc_{M_x}(\xi)$ is homogeneous of degree zero.
For $u \in S$ we have
\begin{align*}
   \Qc_{M_x}(\xi)(u\frac{m}{x^n})&= \Qc_{M_x}(\xi)(\frac{um}{x^n}) \\
  &= \Pc(\xi)(x)^{-n}\Qc_M(\xi)(um)  \\
  &= \Pc(\xi)(x)^{-n}\Pc(\xi)(u)\Qc_M(\xi)(m) \\
  &= \Pc(\xi)(u)\Qc_{M_x}(\xi)(\frac{m}{x^n}).
\end{align*}
Thus $\Qc_{M_x}(\xi)$ is a Steenrod map. The assertion follows from Theorem \ref{char-cartan}.
\end{proof}
Next we give
\begin{proof}[Proof of Theorem \ref{maps}]
We have nothing to show if $M_x = 0$ or $M_{xy} = 0$. So assume that both are non-zero. So $x, xy$ (and hence $y$) are non-nilpotent elements of $S$.
We let $i \colon M_x \rt M_{xy}$ denote the natural map $i(m/x^n) = y^nm/(xy)^n$. This induces a natural map $\eta \colon M_x[[\xi]]^* \rt M_{xy}[[\xi]]^*$.
We have
\begin{align*}
  \Qc_{M_{xy}}(\xi)\circ i(m/x^n) &=  \Qc_{M_{xy}}(\xi)(y^nm/x^ny^n) \\
  &\Pc(\xi)(xy)^{-n}\Qc_M(\xi)(y^nm) \\
   &= \Pc(\xi)(x)^{-n}\Pc(\xi)(y)^{-n}\Pc(\xi)(y)^{n}\Qc_M(\xi)(m) \\
   &=  \Pc(\xi)(x)^{-n}\Qc_M(\xi)(m) \\
   &= \eta\circ\Qc_{M_x}(\xi)(m/x^n).
\end{align*}
From this it easily follows that the map $M_x \rt M_{xy}$ is a Cartan morphism.
\end{proof}
\section{Proof of Theorem \ref{jim}}
We now give
\begin{proof}[Proof of Theorem \ref{jim}]
Let $\mathcal{C}_M = 0 \rt C^0 \rt C^1 \rt \cdots \rt C^m \rt 0$ be the modified \v{C}ech complex on $\xb$ with coefficients in $M$, see \cite[5.1.5]{bs}. Here
$$ C^t = \bigoplus_{1 \leq i_1< i_2 < \cdots<i_t \leq m}M_{x_{i_1}x_{i_2}\cdots x_{i_t}},$$
where the differential $\partial^t \colon C^t \rt C^{t+1}$ is given by \\  $(-1)^{s-1}.\text{nat}\colon M_{x_{i_1}x_{i_2}\cdots x_{i_t}} \rt (M_{x_{i_1}x_{i_2}\cdots x_{i_t}})_{x_{j_s}}$ if
$\{ i_1, \ldots, i_t \} = \{ j_1, \ldots, \wh{j_s}, \ldots, j_{t+1} \}$ and $0$ otherwise. It follows from \ref{loc} and \ref{add} that $C^t$ is a Cartan $S$-module for all $t$. By \ref{maps} and \ref{map-add} it follows that $\partial^t$ is a morphism of Cartan modules. By \ref{complex-cart} it follows that $H^i_I(M)$ has a structure of a Cartan $S$-module.
\end{proof}

\begin{remark}
Suppose $\xb$ and $\yb$ be two different sets of homogeneous generators of $I$. Then  we can give a Cartan $S$-module structure on $H^i_I(M)$ which depends on $\xb$ and another Cartan $S$-module structure which depends on $\yb$. We \emph{do not know} whether these structures are the same. For our applications it suffices to know there is at least one Cartan $S$-module structure on local cohomology modules.
\end{remark}

\begin{remark}
  We note that in our proof of Theorem \ref{jim} the Cartan module $M$ need not be finitely generated as an $S$-module. As a corollary we get that for any homogeneous ideals $I_1, \ldots, I_t$
  the iterated local cohomology module
  $$H^{i_1}_{I_1}(H^{i_2}_{I_2}(\cdots (H^{i_t}_{I_t}(M)) \cdots))$$
  has a structure of a Cartan $S$-module.
\end{remark}
\begin{remark}
In the study of local cohomology the so called Lyubeznik functors plays an important role, see \cite[p.\ 44]{Lyu}. Iterated local cohomology modules forms a large subclass of Lyubeznik functors. We do not know whether if $\mathcal{T}$ is a graded Lyubeznik functor on $\text{Mod}(S)$ then $\mathcal{T}(S)$ is a Cartan module. We however suspect that it is a Cartan $S$-module.
\end{remark}
\s Let $E$ be the $*$-injective hull of $\F_q$ as an $S$-module. Set $(-)^\vee = \ ^*\Hom_S(-, E)$. If $U  $ is  a finitely generated graded $S$-module then $U^\vee$ is a $^*$-Artinian $S$-module.
 If $U  $ is  a  $^*$-Artinian $S$-module then $U^\vee$ is a  finitely generated graded $S$-module. Furthermore the natural map $U \rt U^{\vee \vee}$ is an isomorphism if either $U$ is finitely generated graded or when $U  $ is  a  $^*$-Artinian $S$-module.

\s\label{prev} If $M$ is a finitely generated graded $S$-module then $H^i_{S_+}(M)$ are $^*$-Artinian $S$-module for all $i \geq 0$. Furthermore it is well known that $\dim H^i_{S_+}(M)^\vee$ is $\leq i$, for instance see \cite[10.1]{P}.
\begin{corollary}\label{usha}
  Let $M$ be a finitely generated graded Cartan $S$-module. Set $J_ i = \ann H^i_{S_+}(M)$. Then $\db_{d,0}^{q^t}, \ldots, \db_{d, d-i-1}^{q^t} \in \sqrt{J_i}$ for $i < d$.
\end{corollary}
\begin{proof}
  By \ref{jim} we get that $H^i_{S_+}(M)$ is a Cartan-module. So $J_i$ is a $\Pc^*$-invariant ideal of $S$. Therefore $\sqrt{J_i}$ is also a $\Pc^*$-invariant ideal of $S$, (see \ref{cleo}). Thus $J_i =\ann H^i_{S_+}(M)^\vee$ and the later module has dimension $\leq i$, see \ref{prev}. The result follows from \ref{bella}
\end{proof}
\section{Proof of Theorem \ref{bz-gen}}
We make the convention that the depth of the zero module i $+ \infty$.
In this section we give
\begin{proof}[Proof of Theorem \ref{bz-gen}]
The assertion (ii) $\implies $ (i) is trivial. We prove the converse.
We note that $\db_{d,0}^{q^t},\ldots, \db_{d, d-1}^{q^t}$ is a system of parameters of $S$. So if $M$ is \CM \ then
we have nothing to prove.

Assume $\depth M = d-g$ where $g > 0$. Let $D = D^*(d)^{q^t}$ be the (fractal) Dickson algebra. We consider $M$ as a $D$-module. Set $J =  (\db_{d,d-1}^{q^t}, \cdots, \db_{d,g}^{q^t})$. We note that it suffices to prove $\depth M_P \geq d -g$ for all primes $P$ in $D$ which contain $J$, see \cite[1.2.10]{BH}. The assertion holds if $P = D_+$ by our assumption. Now let $P \supseteq J$ and $P \neq D_+$. Let $D(u)$ be the $*$-canonical module of $D$.
We note that $\dim H^{i}_{S_+}(M)^\vee = \Ext^{d-i}_D(M, D(u))$ has dimension $\leq i$
 (see \ref{usha}).  Set $K_{d-i} = \ann_D \Ext^{d-i}_D(S, D(u))$. So
 $\db_{d,0}^{q^t}, \ldots, \db_{d, d-i-1}^{q^t} \in \sqrt{K_{d-i}}$ by \ref{usha} for $i < d$.

  Let $P \supseteq J$ be minimal. Then $P = J$. As $\db_{d, 0} \notin P$ we get that\\
  $\Ext^{i}_D(M, D(u))_P = 0$ for $i > 0$. So $M_P$ is either \CM \ of $\depth  = d -g$ or $M_P$ is  zero.

  Let $P \supseteq J$ be such that $\height P = d - g + r$ with $r \geq 1$ and that $P \neq D_+$. We assert that $\depth M_P \geq d- g$. If not then $\Ext^{d-g + r -(d -g -i)}_D(M, D(u))_P \neq 0$ for some $i > 0$. So we have $\Ext^{r+i}(M, D(u))_P \neq 0$. We have $U = \Ext^{r+i}(M, D(u)) = H^{d - r - i}_{S_+}(M)^\vee$. So $\db_{d,0}^{q^t}, \cdots, \db_{d, r + i -1}^{q^t} \in \sqrt{\ann_D U}$ by  \ref{usha}. As $U_P \neq 0$ we get that $\db_{d,0}^{q^t}, \cdots, \db_{d, r + i -1}^{q^t} \in P$. If $t+i -1 \geq g$ then note that $D_+ \subseteq P$ which is a contradiction. So  $r + i -1 < g$.
  Thus  $P$ contains $d - g + r + i$ variables where $i \geq 1$. This is a contradiction as $\height P = d - g + r$. The result follows.
\end{proof}
\section{Proof of Theorem \ref{ass}}
\s \label{dil} In this section we give a proof of Theorem \ref{ass}. Throughout $S$ is an unstable $\F_q$-Steenrod algebra and $M$ is a Cartan $S$-module. The following result follows from Theorem \ref{jim}. We give a direct proof.
\begin{lemma}
\label{water}(with hypotheses as in \ref{dil}). Let $I$ be a homogeneous ideal in $S$. Then $H^0_I(M)$ is a Cartan $S$-submodule of $M$.
\end{lemma}
\begin{proof}
  Let $m \in H^0_I(M)$. We show that $\Qc_M^i(m) \in H^0_I(M)$ for all $ i \geq 0$. Say $I^sm = 0$. Set $J = I^s$.
  
   We prove by induction on $i \geq 0$ that
  $u_1u_2\cdots u_{i+1}\Qc^i_M(m) =0$ for any $u_j \in J$. 
  When $i = 0$ we have $\Qc^0_M(m) =m$. So $u_1m = 0$ by our hypothesis. 
  Now assume the result for $i = 0, \ldots, r- 1$ where $r \geq 1$.
  We have $u_1m = 0$. So $\Qc_M^r(u_1m) = 0$. By Cartan identity we have
  $$ 0 = \Qc^r_M(u_1m) = u_1\Qc^r_M(m) + \sum_{i+j=r, j<r}\Pc^i(u_1)\Qc^j_M(m).$$
  Multiply the above equation by $u_2\cdots u_{r+1}$ to obtain
  $u_1u_2\cdots u_{r+1}\Qc^r_M(m) =0$. This will complete the inductive step.
  
  Note we get $I^{(i+1)s}\Qc^i_M(m) = 0$ for all $i \geq 0$. Hence $\Qc_M^i(m) \in H^0_I(M)$.
\end{proof}
We need the following version of a result in \cite[3.13]{Eis}.
\begin{proposition}
  \label{yesu-e}Let $R$ be a Noetherian ring and let $M$ be a finitely generated $R$-module. Let $P$ be a maximal element in $\Ass_R M$ (with respect to inclusion). Then
  $\Ass_R H^0_P(M) = \{ P \}$ and $\Ass_R M/H^0_P(M)  = \Ass_R M \setminus \{ P \}$.
\end{proposition}

We prove Theorem \ref{ass} by induction on $|\Ass_S M|$ (here $|W|$ denotes the cardinality of a set $W$).
The following result enables our induction.
\begin{proposition}
  \label{utah}(with setup as in \ref{dil}). Further assume $M$ is a finitely generated $S$-module. Let $P$ is a maximal element in $\Ass_S M$ (with respect to inclusion). Then $P$ is a $\Pc^*$-invariant prime ideal of $S$.
\end{proposition}
\begin{proof}
  Set $N = H^0_P(M)$. Then $N$ is a Cartan submodule of $M$, see \ref{water}. So $J = \ann_S N$ is a $\Pc^*$ invariant ideal of $S$, see \ref{ann-p*}. So $\sqrt{J}$ is also a $\Pc^*$ invariant ideal of $S$, see \ref{cleo}. 
  
  As $N$ is a finitely generated $S$-module and $P$-torsion we get that $P^sN = 0$ for some $s \geq 1$. Thus $P^s \subseteq J$.
  As $P \in \Ass_S N$ we get that $P = (0 \colon n) $ for some $n \in N$. We have $Jn = 0$. So $J \subseteq P$.
  We have $P^s \subseteq J \subseteq P$. So $P = \sqrt{J}$ and hence a $\Pc^*$-invariant ideal of $S$.
\end{proof}
We now give
\begin{proof}[Proof of Theorem \ref{ass}] We prove the result by induction on $|\Ass_S M|$.

We first consider the case when $|\Ass_S M| = 1$. Let $\{ P \} = \Ass_S M$. Then trivially $P$ is a maximal element of $\Ass_S M$ and so by \ref{utah} we get that $P$ is $\Pc^*$-invariant.

Now consider the case when $|\Ass_S M| = r$ and the result is known for all finitely generated Cartan $S$-modules $N$ with $|\Ass_S N| = r -1$.
Let $P$ be a maximal element of $\Ass_S M$. Then by \ref{utah} we get $P$ is a $\Pc^*$-invariant prime ideal of $S$. By \ref{water} we get that $H^0_P(M)$ is a Cartan submodule of $M$. So $N = M/H^0_P(M)$ is also a finitely generated Cartan $S$-module. By \ref{yesu-e} we get that $\Ass_S N = \Ass_S M \setminus \{P \}$. So $|\Ass_S N| = r -1$ and by our induction hypotheses we get that $Q$ is a $\Pc^*$-invariant prime ideal of $S$ for every associate prime $Q$ of $N$. Our assertion follows.
\end{proof}
We give an example which shows the result in Theorem \ref{ass} need not hold if $M$ is not a finitely generated Cartan $S$-module.
\begin{example}\label{non-ass}
  Let $R = \F_q[X_1, \ldots, X_d]$ with $d \geq 2$ and the standard Steenrod algebra structure. As discussed earlier $R$ has infinitely many  homogeneous prime ideals which are not $\Pc^*$-invariant. Choose one such prime $P$ with height $g \geq 1$. Then $M = H^g_P(S)$ is a not finitely generated $R$-module with $\Ass M = \{ P\}$. By \ref{jim}
  we get that $M$ is a Cartan $S$-module.
\end{example}
\section{Proof of Applications III}
\s\label{mad} In this section we give all the proofs of result stated in our Application III. Throughout $S$ is a $\F_q$-Steenrod algebra and $M$ is a finitely generated Cartan $S$-module. 

We will need the following result proved in \cite[8.2]{P} for invariant rings. The same proof works in general.
\begin{lemma}
\label{eclair} Let $M$ be a Cartan $S$-module. Assume $P$ is not a $\Pc^*$-invariant prime ideal. If $P^r \subseteq \ann M$ then $M_P = 0$.
\end{lemma}
\
The following result follows easily from Theorem \ref{jim} and \ref{eclair}.
\begin{corollary}
\label{stock}(with setup as in \ref{mad}). Let $P$ be a prime ideal in $S$ which is \emph{not} $\Pc^*$-invariant. If $H^i_P(M)$ is a finitely generated module as an $S$-module then $H^i_P(M)_P = 0$. \qed
\end{corollary}
The following result is crucial for our applications.
\begin{theorem}\label{crucial}
(with setup as in \ref{mad}). Assume further that $S$ is equi-dimensional and $\Supp(M) = \Spec(S)$. Let $r \geq 0$. The following assertions are equivalent:
\begin{enumerate}[\rm (i)]
\item
$\depth M_P \geq \min \{ \height P, r \}$ for all homogeneous primes $P$ of $S$.
\item
$\depth M_P \geq \min \{ \height P, r \}$ for all homogeneous $\Pc^*$-invariant primes $P$ of $S$.
\end{enumerate}
\end{theorem}
\begin{proof}
The assertion (i) $\implies$ (ii) is trivial.
We now prove (ii) $\implies$ (i). We have nothing to show when $r = 0$. So assume $r \geq 1$.

 Let $P$ be a homogeneous  prime ideal of $S$.
Assume that $\height P  = t \geq 1$. Set
$$f_P(M) = \min \{ i \mid H^i_P(M) \ \text{is not finitely generated $S$-module} \}. $$
Set $s_P = \min \{ \height P, r \}$.
We make the following

Claim: If $P$ is  a homogeneous prime ideal then $f_P(M)\geq  s_P$.

We first prove that the Claim proves the result. Let  $P$ be not a $\Pc^*$-invariant homogeneous  ideal. Then $H^i_P(M)$ are finitely generated $S$-modules for $i < s_P$.  By \ref{stock} it follows that $H^i_{PS_P}(M_P) = H^i_{P}(M)_P = 0$ for $i < s_P$. Thus $\depth M_P \geq s_P$. The result follows.

We note that the assertion $f_P(M) \geq s_P $ holds trivially when $r = 1$.

Next consider the case when $r \geq 2$ and the result holds for $\leq r - 1$.
By induction hypothesis, $f_P(M) \geq \min \{ \height P, r - 1 \}$ for all homogeneous primes.   In particular $\depth M_\q \geq  \min \{\height \q, r - 1 \}$ for all homogeneous non $\Pc^*$-invariant prime ideals $\q$ of $S$. We also have $\depth  M_Q \geq s_Q$ for all $\Pc^*$-invariant prime ideals $Q$ by our assumption.
We may assume $\height P \geq r$, otherwise there is nothing to show.
 By \cite[13.1.17]{bs} we have
\[
f_P(M) = \min \{ \depth M_\q  + \height((P+\q)/\q) \mid \q \in \  ^* \Spec(S) \setminus V(P) \}.
\]
If $\height \q \leq r -1$ then $\depth M_\q = \height \q$. We have
\begin{align*}
  \depth M_\q + \height((P+\q)/\q) &=  \height \q + \height((P + \q)/P) \\
  &= \height \q + \height(P + \q) - \height \q,  \\
  &= \height(P + \q)\geq  r.
\end{align*}
The second equality holds as $S$ is equi-dimensional and catenary, see  \cite[Lemma 2, p.\ 250]{Ma} for the local version; same proof works for the *-local case.
If $\height \q \geq r$ then $$\depth M_\q  + \height((P+\q)/\q)  \geq r - 1 + 1 = r. $$
Thus $f_P(M) \geq r$. So Claim follows.
\end{proof}

Surprisingly  the $S_r$ property of a Cartan-module only depends on $\Pc^*$-invariant prime ideals in the following sense.
\begin{theorem}\label{sr}(with hypotheses as in \ref{crucial}).  The following conditions are equivalent:
  \begin{enumerate}[\rm (i)]
    \item $M$ satisfies $S_r$ property.
    \item $\depth M_P \geq \min\{\height P, r \}$  for all $\Pc^*$-invariant homogeneous primes $P$ of $S$.
  \end{enumerate}
\end{theorem}
\begin{proof}
The assertion (i) $\implies$ (ii) is clear.

We now prove (ii) $\implies$ (i). It is well-known that it  suffices to prove $\depth M_P \geq \min \{\height P, r \}$ for homogeneous primes in $S$, for instance  see \cite[10.3]{P} for the case when $M = S$.
By \ref{crucial} it follows that $\depth M_P \geq \min \{ \height P, r \}$ for all homogeneous primes $P$ of $S$. The result follows.
\end{proof}

Next we give
\begin{proof}[Proof of Theorem \ref{s_r}]
We first prove (ii) $\implies $ (i). Let $P$ be a graded $\Pc^*$-invariant prime ideal of $S$. We note that if  $\height P = g$ then $\dim S/P = d - g$, as $S$ is equi-dimensional and universally catenary, see \cite[Lemma 2, p.\ 250]{Ma} for the local version; same proof works for the *-local case. By \ref{bella} we get that if $g \leq r$ then $P$ contains $\db_{d,0}^{q^t}, \ldots, \db_{d,g-1}^{q^t}$.
It follows that $M_P$ is \CM. If $\height P \geq r + 1$ then by \ref{bella}  we get that  $P$ contains $\db_{d,0}, \ldots, \db_{d,r-1}$.
So $\depth M_P \geq r$. The result follows from \ref{sr}.

Next we prove  (i) $\implies $ (ii). The extension $D^*(d)^{q^t} \subseteq S$ is a finite extension. As $(\db_{d,0}^{q^t}, \ldots, \db_{d,r-1}^{q^t})$ is a prime ideal of height $r$ in $D^*(d)$, it follows that $(\db_{d,0}^{q^t}, \ldots, \db_{d,r-1}^{qt})S$ has height $r$ as an ideal in $S$. As $M$ satisfies $S_r$ it follows from \cite[1.2.10]{BH} that $\grade \left( (\db_{d,0}, \ldots, \db_{d,r-1})S, M \right)) = r$. The result follows.
\end{proof}

\section{proof of Corollary \ref{cm}}
In this section we give
\begin{proof}[Proof of Corollary \ref{cm}]
By \ref{sbi} we may assume $\dim M = \dim S = d$.
  We work with the setup as in \ref{lannes}.  Set $D = D(d)^{q^t}$. Let $D(u)$ be the canonical module of $D$. By local duality we get that $H^i_{\S_+}(M)^\vee = \Ext^{d-i}_D(M, D(u)) $. So $\db_{d,0}^{q^t} \in \sqrt{\ann\Ext^{d-i}_D(M, D(u))}$ for $i < d$. We note that we have a finite extension $D_{\db_{d,0}^{q^t} }\hookrightarrow S_{\db_{d,0}^{q^t} }$. So it suffices to show that $M_{\db_{d,0}^{q^t}}$ is \CM \ as a $D_{\db_{d,0}^{q^t}}$-module. Let $P$ be any homogeneous prime  in $D$ not containing $\db_{d,0}^{q^t}$. Then it follows that $\Ext^{i}_D(M, D(u))_P = 0$ for $i > 0$. It follows that $M_P$ is \CM. The result follows.
\end{proof}

\section{Proof of Theorem \ref{mono}}
In this section we prove Theorem \ref{mono}.
To prove this result we need Theorem 23.3 from \cite{Ma}.   Unfortunately  there is a typographical error in the statement of Theorem 23.3 in \cite{Ma}.
So we state it here.
\begin{theorem}\label{23}
Let $\varphi \colon A \rt B$ be a homomorphism of Noetherian rings, and let $E$ be an $A$-module and $G$ a $B$-module. Suppose that
$G$ is flat over $A$; then we have the following:
\begin{enumerate}[\rm (i)]
\item
if $\fp \in \Spec A$ and $G/\fp G \neq 0$ then
\[
 ^a \varphi \left( \Ass_B(G/\fp G)  \right) = \Ass_A (G/\fp G) = \{\fp \}.
\]
\item
$\displaystyle{\Ass_B(E\otimes_A G) = \bigcup_{\fp \in \Ass_A(E) } \Ass_B(G/\fp G).}$
\end{enumerate}
\end{theorem}
\begin{remark}
In  \cite{Ma} $\Ass_A(E\otimes G)$ is typed instead of $\Ass_B(E\otimes G)$. Also note that
$^a \varphi(\fP) = \fP \cap A$ for $\fP \in \Spec B$.
\end{remark}
\begin{proof}
The inclusion $R \rt T$ is faithfully flat. Also note that as $I$ is a monomial ideal it is extended from $R$. Say $I = LT$ for a monomial ideal $L$ of $R$. Set $N = H^i_J(R/L)$. Then $M = N \otimes_R T$. As $R \rt T$ is faithfully flat we get that $N$ is a finitely generated $R$-module. By \ref{s-mono} we get that $L$ is a $\Pc^*$-invariant ideal of $R$. So $R/L$ is a Cartan $R$-module.  By \ref{jim} we get that $N$ is a Cartan $R$-module. It is also finitely generated as an $R$-module. Thus if $P \in \Ass_R N$ then $P$ is a $\Pc^*$-invariant ideal in $R$. So by \cite[11.4.1]{SL} we get that $P$ is generated by some linear homogeneous polynomials of $R$. In \ref{23} we take $A = R, B = T$ and $G = T$. We note that if $P \in \Spec(R)$ is homogeneous and generated by linear polynomials then $PT$ is a prime in $T$. The result follows from \ref{23}.
\end{proof}

\section{Proof of Theorem \ref{bound}}
In this section we give
\begin{proof}[Proof of Theorem \ref{bound}]
Set $t = \depth S_P$.
Suppose if possible $f_P > t$. Then note that $H^t_P(S)$ is finitely generated $S$-module. Then by \ref{stock} we get that $H^t_P(S)_P = 0$.
But $H^t_{PS_P}(S_P) = H^t_P(S)_P = 0$. This is a  contradiction as $t = \depth S_P$.
\end{proof}
We give an example which shows that the bound in Theorem \ref{bound} need not hold for $\Pc^*$-invariant ideals.
\begin{example}\label{counter-bound}
  Let $I$ be a monomial ideal in $R = \F_\q[X_1, \ldots, X_n]$ such that $S = R/I$ is generalized \CM \ but not \CM. We note that $S$ is an unstable Steenrod 
  $\F_q$-algebra, see \ref{s-mono}. We note that $f_{S_+} = d = \dim S$ while $\depth S_{S_+} < \dim S_{S_+} = d$. For a specific example see \cite[3.5]{TY}.
\end{example}


\begin{thebibliography} {99}



\bibitem{bz}
D.~Bourguiba and S.~Zarati,
\emph{Depth and the Steenrod algebra},
With an appendix by J.~Lannes.
Invent. Math. 128 (1997), no. 3, 589--602.

\bibitem{bs}  M. Brodmann and R. Y. Sharp,
Local Cohomology: An algebraic introduction with geometric applications,
Cambridge University Press, 60, 1998.

\bibitem{BH}  W. Bruns and J. Herzog,
Cohen-Macaulay Rings, revised edition,
Cambridge Studies in Advanced Mathematics, 39.
Cambridge University Press, 1998.

\bibitem{Eis}
D.~Eisenbud, 
\emph{Commutative algebra,
With a view toward algebraic geometry},
Grad. Texts in Math., 150
Springer-Verlag, New York, 1995.

\bibitem{LS}
P.~S.~Landweber and R.~E.~Stong,
\emph{The depth of rings of invariants over finite fields},
Number theory (New York, 1984–1985), 259--274.
Lecture Notes in Math., 1240
Springer-Verlag, Berlin, 1987.

\bibitem{Lyu}
G.~Lyubeznik,
 \emph{Finiteness properties of local cohomology modules (an application of  $D$ -modules to commutative algebra)},
Invent. Math. 113 (1993), no. 1, 41--55.


\bibitem{Ma}
H.~Matsumura, \emph{Commutative ring theory}, Cambridge Studies in Advanced
  Mathematics, vol.~8, Cambridge University Press, Cambridge, 1986.


\bibitem{NS}
M.~D.~Neusel and  L.~Smith,
\emph{Invariant theory of finite groups},
Math. Surveys Monogr., 94
American Mathematical Society, Providence, RI, 2002.

\bibitem{P}
T.~J.~Puthenpurakal,
\emph{Annhilators of local cohomology modules over modular invariant rings and Dickson polynomials}.
Preprint,	arXiv:2603.16491.

\bibitem{SL}
L.~Smith,
\emph{Polynomial invariants of finite groups},
Res. Notes Math., 6
A K Peters, Ltd., Wellesley, MA, 1995.

\bibitem{TY}
N.~Terai and K.~Yoshida, 
\emph{Buchsbaum Stanley-Reisner rings with minimal multiplicity},
Proc. Amer. Math. Soc. 134 (2006), no. 1, 55--65.



\end{thebibliography}
\end{document}